\documentclass[a4paper,runningheads]{llncs}
\usepackage[latin9]{inputenc}
\usepackage{float}
\usepackage{amsmath}
\usepackage{amssymb}
\usepackage{graphicx}

\makeatletter

\floatstyle{ruled}
\newfloat{algorithm}{tbp}{loa}
\providecommand{\algorithmname}{Algorithm}
\floatname{algorithm}{\protect\algorithmname}

\@ifundefined{date}{}{\date{}}

\pdfoutput=1
\usepackage{url}\urldef{\mailsa}\path|{alfred.hofmann, ursula.barth, ingrid.haas, frank.holzwarth,|
\urldef{\mailsb}\path|anna.kramer, leonie.kunz, christine.reiss, nicole.sator,|
\urldef{\mailsc}\path|erika.siebert-cole, peter.strasser, lncs}@springer.com|    
\newcommand{\keywords}[1]{\par\addvspace\baselineskip
\noindent\keywordname\enspace\ignorespaces#1}

\makeatother

\begin{document}
\mainmatter 


\title{Linear Superiorization for Infeasible Linear Programming}

\author{Yair Censor and Yehuda Zur}

\authorrunning{Yair Censor and Yehuda Zur}

\institute{Department of Mathematics, University of Haifa, \linebreak{}
 Mt. Carmel, Haifa 3498838, Israel}

\title{Linear Superiorization for Infeasible Linear Programming}

\maketitle




\toctitle{Lecture Notes in Computer Science}

\tocauthor{Authors' Instructions}
\begin{abstract}
Linear superiorization (abbreviated: LinSup) considers linear programming
(LP) problems wherein the constraints as well as the objective function
are linear. It allows to steer the iterates of a feasibility-seeking
iterative process toward feasible points that have lower (not necessarily
minimal) values of the objective function than points that would have
been reached by the same feasiblity-seeking iterative process without
superiorization. Using a feasibility-seeking iterative process that
converges even if the linear feasible set is empty, LinSup generates
an iterative sequence that converges to a point that minimizes a proximity
function which measures the linear constraints violation. In addition,
due to LinSup's repeated objective function reduction steps such a
point will most probably have a reduced objective function value.
We present an exploratory experimental result that illustrates the
behavior of LinSup on an infeasible LP problem. \keywords{Superiorization;
perturbation resilience; infeasible linear programming; feasibility-seeking;
simultaneous projection algorithm; Cimmino method; proximity function.} 
\end{abstract}

\section{Introduction: The General Concept of Superiorization}

Given an algorithmic operator $\mathcal{A}:X\rightarrow X$ on a Hilbert
space $X$, consider the iterative process

\begin{equation}
x^{0}\in X,\;x^{k+1}=\mathcal{A}\left(x^{k}\right),\;\mathrm{for}\,\mathrm{all}\;k\geqslant0,\label{eq:iter-proc-basic}
\end{equation}
and let $SOL\left(P\right)$ denote the solution set of some problem
$P$ of any kind. The iterative process is said to solve $P$ if,
under some reasonable conditions, any sequence $\left\{ x^{k}\right\} _{k=0}^{\infty}$
generated by the process converges to some $x^{*}\in SOL\left(P\right)$.
An iterative process (\ref{eq:iter-proc-basic}) that solves $P$
is called perturbation resilient if the process

\begin{equation}
y^{0}\in X,\;y^{k+1}=\mathcal{A}\left(y^{k}+v^{k}\right),\;\mathrm{for}\,\mathrm{all}\;k\geqslant0,\label{eq:iter-proc-basic-1}
\end{equation}
also solves $P$, under some reasonable conditions on the sequence
of perturbation vectors $\left\{ v^{k}\right\} _{k=0}^{\infty}\subseteq X$.
The iterative processes of (\ref{eq:iter-proc-basic}) and (\ref{eq:iter-proc-basic-1})
are called ``the basic algorithm'' and ``the superiorized version
of the basic algorithm'', respectively.

Superiorization aims at identifying perturbation resilient iterative
processes that will allow to use the perturbations in order to steer
the iterates of the superiorized algorithm so that, while retaining
the original property of converging to a point in $SOL\left(P\right)$,
they will also do something additional useful for the original problem
$P$, such as converging to a point with reduced values of some given
objective function. These concepts are rigorously defined in several
recent works in the field, we refer the reader to the recent reviews
\cite{h14super},\cite{censor15weak} and references therein. More
material about the current state of superiorization can be found also
in \cite{linsup16}, \cite{hgdc12} and \cite{reem15new}.

A special case of prime importance and significance of the above is
when $P$ is a convex feasibility problem (CFP) of the form: Find
a vector $x^{*}\in\cap_{i=1}^{I}C_{i}$ where $C_{i}\subseteq R^{J}$,
the $J$-dimensional Euclidean space, are closed convex subsets, and
the perturbations in the superiorized version of the basic algorithm
are designed to reduce the value of a given objective function $\phi$.

In this case the basic algorithm (\ref{eq:iter-proc-basic}) can be
any of the wide variety of feasibility-seeking algorithms, see, e.g.,
\cite{bb96}, \cite{cccdh10} and \cite{annotated}, and the perturbations
employ nonascent directions of $\phi$. Much work has been done on
this as can be seen in the Internet bibliography at \cite{bib-page}.

The usefulness of this approach is twofold: First, feasibility-seeking
is, on a logical basis, a less-demanding task than seeking a constrained
minimization point in a feasible set. Therefore, letting efficient
feasibility-seeking algorithms ``lead'' the algorithmic effort and
modifying them with inexpensive add-ons works well in practice.

Second, in some real-world applications the choice of an objective
function is exogenous to the modeling and data acquisition which give
rise to the constraints. Thus, sometimes the limited confidence in
the usefulness of a chosen objective function leads to the recognition
that, from the application-at-hand point of view, there is no need,
neither a justification, to search for an exact constrained minimum.
For obtaining ``good results'', evaluated by how well they serve
the task of the application at hand, it is often enough to find a
feasible point that has reduced (not necessarily minimal) objective
function value\footnote{Some support for this reasoning may be borrowed from the American
scientist and Noble-laureate Herbert Simon who was in favor of ``satisficing''
rather then ``maximizing''. Satisficing is a decision-making strategy
that aims for a satisfactory or adequate result, rather than the optimal
solution. This is because aiming for the optimal solution may necessitate
needless expenditure of time, energy and resources. The term ``satisfice''
was coined by Herbert Simon in 1956 \cite{simon1956satisficing},
see: https://en.wikipedia.org/wiki/Satisficing.}.

\section{Linear Superiorization}

\subsection{The problem and the algorithm}

Let the feasible set $M$ be 
\begin{equation}
M:=\{x\in R^{J}\mid Ax\leq b,\text{ }x\geq0\}\label{eq:linear-feas-prob}
\end{equation}
where the $I\times J$ real matrix $A=(a_{j}^{i})_{i=1,j=1}^{I,J}$
and the vector $b=(b_{i})_{i=1}^{I}\in R^{I}$ are given.

For a basic algorithm we pick a feasibility-seeking projection method.
Here projection methods refer to iterative algorithms that use projections
onto sets while relying on the general principle that when a family
of, usually closed and convex, sets is present, then projections onto
the individual sets are easier to perform than projections onto other
sets (intersections, image sets under some transformation, etc.) that
are derived from the individual sets.

Projection methods may have different algorithmic structures, such
as block-iterative projections (BIP) or string-averaging projections
(SAP) (see, e.g., the review paper \cite{censor-segal2008} and references
therein) of which some are particularly suitable for parallel computing,
and they demonstrate nice convergence properties and/or good initial
behavior patterns.

This class of algorithms has witnessed great progress in recent years
and its member algorithms have been applied with success to many scientific,
technological and mathematical problems. See, e.g., the 1996 review
\cite{bb96}, the recent annotated bibliography of books and reviews
\cite{annotated} and its references, the excellent book \cite{CEG12},
or \cite{cccdh10}.

An important comment is in place here. A CFP can be translated into
an unconstrained minimization of some proximity function that measures
the feasibility violation of points. For example, using a weighted
sum of squares of the Euclidean distances to the sets of the CFP as
a proximity function and applying steepest descent to it results in
a simultaneous projections method for the CFP of the Cimmino type.
However, there is no proximity function that would yield the sequential
projections method of the Kaczmarz type, for CFPs, see \cite{baillon2012}.

Therefore, the study of feasibility-seeking algorithms for the CFP
has developed independently of minimization methods and it still vigorously
does, see the references mentioned above. Over the years researchers
have tried to harness projection methods for the convex feasibility
problem to LP in more than one way, see, e.g., Chinneck's book \cite{C08}.

The mini-review of relations between linear programming and feasibility-seeking
algorithms in \cite[Section 1]{nurminski2015single} sheds more light
on this. Our work in \cite{linsup16} and here leads us to study whether
LinSup can be useful for either feasible or infeasible LP problems.

The objective function for linear superiorization will be 
\begin{equation}
\phi(x):=\left\langle c,x\right\rangle \label{eq:lin-objective}
\end{equation}
where $\left\langle c,x\right\rangle $ is the inner product of $x$
and a given $c\in R^{J}.$

In the footsteps of the general principles of the superiorization
methodology, as presented for general objective functions $\phi$
in previous publications, we use the following linear superiorization
(LinSup) algorithm. The algorithm and its implementation details follow
closely those of \cite{linsup16} wherein only feasible constraints
were discussed.

The input to the algorithm consists of the problem data $A,$ $b,$
and $c$ of (\ref{eq:linear-feas-prob}) and (\ref{eq:lin-objective}),
respectively, a user-chosen initialization point $\bar{y}$ and a
user-chosen parameter (called here kernel) $0<\alpha<1$ with which
the algorithm generates the step-sizes $\beta_{k,n}$ by the powers
of the kernel $\eta_{\ell}=\alpha^{\ell},$ as well as an integer
$N$ that determines the quantity of objective function reduction
perturbation steps done per each feasibility-seeking iterative sweep
through all linear constraints. The perturbation direction $-\frac{{\displaystyle c}}{{\displaystyle \left\Vert c\right\Vert _{2}}}$
used in step 10 of Algorithm 1 is a nonascend direction of the linear
objective function, as required by the general principles of the superiorization
methodology, see, e.g., \cite[Subsection II.D]{hgdc12}.

\begin{algorithm}[H]
\label{alg_super}\textbf{Algorithm 1. The Linear Superiorization
(LinSup) Algorithm} 
\end{algorithm}

\begin{enumerate}
\item \textbf{set} $k=0$ 
\item \textbf{set} $y^{k}=\bar{y}$ 
\item \textbf{set }$\ell_{-1}=0$ 
\item \textbf{while }stopping rule not met\textbf{ do} 
\item $\qquad$\textbf{set} $n=0$ 
\item \ \ \ \ \ \ \textbf{set} $\ell=rand(k,\ell_{k-1})$ 
\item $\qquad$\textbf{set} $y^{k,n}=y^{k}$ 
\item $\qquad$\textbf{while }$n$\textbf{$<$}$N$ \textbf{do} 
\item $\qquad\qquad$\textbf{set} $\beta_{k,n}=\eta_{\ell}$ 
\item $\qquad\qquad$\textbf{set} $z=y^{k,n}-\beta_{k,n}\frac{{\displaystyle c}}{{\displaystyle \left\Vert c\right\Vert _{2}}}$ 
\item $\qquad\qquad$\textbf{set }$n\leftarrow n+1$ 
\item $\qquad\qquad$\textbf{set }$y^{k,n}$\textbf{$=$}$z$ 
\item \ \ \ \ \ \ \ \ \ \ \ \ \textbf{set }$\ell\leftarrow\ell+1$ 
\item \ \ \ \ \ \ \textbf{end while} 
\item \ \ \ \ \ \ \textbf{set }$\ell_{k}=\ell$ 
\item \qquad{}\textbf{set }$y^{k+1}$\textbf{$=$}$\mathcal{A}\left(y^{k,N}\right)$ 
\item \qquad{}\textbf{set }$k\leftarrow k+1$ 
\item \textbf{end while} 
\end{enumerate}
All quantities in this algorithm are detailed and explained below,
except for the choice of the basic algorithm for the feasibility-seeking
operator represented by $\mathcal{A}$ in step 16 of Algorithm \ref{alg_super}
which appear in the next subsection.

\textbf{Step-sizes of the perturbations}.\textbf{ }The step sizes
$\beta_{k,n}$ in Algorithm 1 must be such that $0<\beta_{k,n}\leq1$
in a way that guarantees that they form a summable sequence $\sum_{k=0}^{\infty}\sum_{n=0}^{N-1}\beta_{k,n}<\infty,$
see, e.g., \cite{censor-zas2015}. To this end Algorithm 1 assumes
that we have available a summable sequence $\{\eta_{\ell}\}_{\ell=0}^{\infty}$
of positive real numbers generated by $\eta_{\ell}=\alpha^{\ell}$
, where $0<\alpha<1$. Simultaneously with generating the iterative
sequence $\{y^{k}\}_{k=0}^{\infty},$ a subsequence of\ $\{\eta_{\ell}\}_{\ell=0}^{\infty}$
is used to generate the step sizes $\beta_{k,n}$ in step 9 of Algorithm
1. The number $\alpha$ is called the kernel of the sequence $\{\eta_{\ell}\}_{\ell=0}^{\infty}.$

\textbf{Controlling the decrease of the step-sizes of objective function
reduction}. If during the application of Algorithm 1 the step sizes
$\beta_{k,n}$ decrease too fast then too little leverage is allocated
to the objective function reduction activity that is interlaced into
the feasibility-seeking activity of the basic algorithm. This delicate
balance can be controlled by the choice of the index $\ell$ updates
and separately by the value of $\alpha$ whose powers $\alpha^{\ell}$
determine the step sizes $\beta_{k,n}$ in step 9. In our work we
adopt a strategy for updating the index $\ell$ that was proposed
and implemented for total variation (TV) image reconstruction from
projections by Prommegger and by Langthaler in \cite[page 38 and Table 7.1 on page 49]{prommegger2014}
and in \cite{langthaler2014}, respectively. This strategy advocates
to set $\ell$ at the beginning of every new iteration sweep (steps
5 and 6) to a random number between the current iteration index $k$
and the value of $\ell$ from the last iteration sweep, i.e., $\ell_{k}=rand(k,\ell_{k-1}).$

\textbf{The proximity function}. To measure the feasibility-violation
(or level of disagreement) of a point with respect to the target set
$M$ we used the following proximity function 
\begin{equation}
\Pr(x):=\frac{1}{2I}{\displaystyle \sum_{i=1}^{I}}\frac{\left(\left(\left\langle a^{i},x\right\rangle -b_{i}\right)_{+}\right)^{2}}{{\displaystyle \sum\limits _{j=1}^{J}}\left(a_{j}^{i}\right)^{2}}+\frac{1}{2J}{\displaystyle \sum\limits _{j=1}^{J}}\left(\left(-x_{j}\right)_{+}\right)^{2}\label{eq:prox}
\end{equation}
where the plus notation means, for any real number $d,$ that $d_{+}:=\max(d,0).$

\textbf{The number }$N$\textbf{ of perturbation steps}. This number
$N$ of perturbation steps that are performed prior to each application
of the feasibility-seeking operator $\mathcal{A}$ (in step 16) affects
the performance of the LinSup algorithm. It influences the balance
between the amounts of computations allocated to feasibility-seeking
and those allocated to objective function reduction steps. A too large
$N$ will make Algorithm 1 spend too much resources on the perturbations
that yield objective function reduction.

\textbf{Handling the nonnegativity constraints}. The nonnegativity
constraints in (\ref{eq:linear-feas-prob}) are handled by projections
onto the nonnegative orthant, i.e., by taking the iteration vector
in hand after each iteration of Cimmino's feasibility-seeking algorithm
applied to all $I$ row-inequalities of (\ref{eq:linear-feas-prob})
and setting its negative components to zero while keeping the others
unchanged.

\subsection{Cimmino's feasibility-seeking algorithm as the basic algorithm\label{subsec:AMS}}

We use the simultaneous projections method of Cimmino for linear inequalities,
see, e.g. \cite{LAA1985cimmino}, as the basic algorithm for the feasibility-seeking
operator represented by$\mathcal{A}$ in step 16 of Algorithm 1. Denoting
the half-spaces represented by individual rows of (\ref{eq:linear-feas-prob})
by $H_{i},$ 
\begin{equation}
H_{i}:=\{x\in R^{J}\mid\left\langle a^{i},x\right\rangle \leq b_{i}\},
\end{equation}
where $a^{i}\in R^{J}$ is the $i$-th row of $A$ and $b_{i}\in R$
is the $i$-th component of $b$ in (\ref{eq:linear-feas-prob}),
he orthogonal projection of an arbitrary point $z\in R^{J}$ onto
$H_{i},$ has the closed-form 
\begin{equation}
P_{H_{i}}(z)=\left\{ \begin{array}{ll}
z-{\displaystyle \frac{\left\langle a^{i},z\right\rangle -b_{i}}{\Vert a^{i}\Vert^{2}}a^{i},} & \text{if }\left\langle a^{i},z\right\rangle >b_{i},\\
z, & \text{if }\left\langle a^{i},z\right\rangle \leq b_{i}.
\end{array}\right.
\end{equation}

\newpage{}

\begin{algorithm}[H]
\label{alg:cimmino}\textbf{Algorithm 2. The Simultaneous Feasibility-Seeking
Projection Method of Cimmino} 
\end{algorithm}

\textbf{Initialization}: $x^{0}\in R^{J}$ is arbitrary.

\textbf{Iterative step}: Given the current iteration vector $x^{k}$
the next iterate is calculated by 
\begin{equation}
x^{k+1}=x^{k}+\lambda_{k}\left(\sum_{i=1}^{I}w_{i}\left(P_{H_{i}}(x^{k})-x^{k}\right)\right)
\end{equation}

with weights $w_{i}\geq0$ for all $i\in I,$ and $\sum_{i=1}^{I}w_{i}=1.$

\textbf{Relaxation parameters}: The parameters $\lambda_{k}$ are
such that $\epsilon_{1}\leq\lambda_{k}\leq2-\epsilon_{2},$ for all
$k\;\geq\;0,$ with some, arbitrarily small, fixed, $\epsilon_{1},\epsilon_{2}>0.$

This Cimmino simultaneous feasibility-seeking projection algorithm
is known to generate convergent iterative sequences even if the intersection
$\cap_{i=1}^{I}H_{i}$ is empty, as the following, slightly paraphrased,
theorem tells. 
\begin{theorem}
\cite[Theorem 3]{LAA1985cimmino} For any starting point $x^{0}\in R^{J}$,
any sequence $\{x^{k}\}_{k=0}^{\infty},$ generated by the simultaneous
feasibility-seeking projection method of Cimmino (Algorithm 2) converges.
If the underlying system of linear inequalities is consistent, the
limit point is a feasible point for it. Otherwise, the limit point
minimizes $f(x):=\sum_{i=1}^{I}w_{i}\parallel P(x)-x\parallel^{2}$,
i.e., it is a weighted (with the weights $w_{i}$) least squares solution
of the system. 
\end{theorem}

\section{An Empirical Result}

Employing MATLAB 2014b \cite{matlab}, we created five test problems
each with 2500 linear inequalities in $R^{J},\:J=2000$. The entries
in 1250 rows of the matrix $A$ in (\ref{eq:linear-feas-prob}) were
uniformly distributed random numbers from the interval $(-1,1).$
The remaining 1250 rows were defined as the negatives of the first
1250 rows, i.e., $a_{j}^{1250+t}=-a_{j}^{t}$ for all $t=1,2,\ldots,1250$
and all $j=1,2,\ldots,2000.$ This guarantees that the two sets of
rows represent parallel half-spaces with opposing normals. For the
right-hand side vectors, the components of $b$ associated with the
first set of 1250 rows in (\ref{eq:linear-feas-prob}) were uniformly
distributed random numbers from the interval $(0,100).$ The remaining
1250 components of each $b$ were chosen as follows: $b_{1250+t}=-b_{t}-rand(100,200)$
for all $t=1,2,\ldots,1250$. This guarantees that the distance between
opposing parallel half-spaces is large making them inconsistent, i.e.,
having no point in common, and that the whole system is infeasible.

For the linear objective function, the components of $c$ were uniformly
distributed random numbers from the interval $(-2,1).$ All runs of
Algorithm 1 and Algorithm 2 were initialized at $\bar{y}=10\cdot\mathbf{1}$
and $x^{0}=10\cdot\mathbf{1}$, respectively, where $\mathbf{1}$
is the vector of all 1's.

We ran Algorithm 1 on each problem until it ceased to make progress,
by using the stopping rule 
\begin{equation}
\frac{{\displaystyle \left\Vert y^{k}-y^{k-1}\right\Vert }}{\left\Vert y^{k}\right\Vert }\leq10^{-4}.
\end{equation}
The same stopping rule was used for runs of Algorithm 2. The relaxation
parameters in Cimmino's feasibility-seeking basic algorithm in step
16 of Algorithm 1 were fixed with $\lambda_{k}=1.99$ for all $k\geqslant0.$
Based on our work in \cite{linsup16} we used $N=20$ and $\alpha=0.99$
in steps 8 and 9 of Algorithm 1, respectively, where $\eta_{\ell}=\alpha^{\ell}.$

The three figures, presented below, show results for the five different
(but similarly generated) families of inconsistent linear inequalities
along with nonnegativity constraints. Figures 1 and 2, in particular,
show that the perturbation steps 5-15 of the LinSup Algorithm 1 initially
work and reduce the objective function value powerfully during the
first ca. 500 iterative sweeps (an iterative sweep consists of one
pass through steps 5-17 in Algorithm 1 or one pass through all linear
inequalities and the nonnegativity constraints in Algorithm 2). As
iterative sweeps proceed the perturbations in Algorithm 1 loose steam
because of the decreasing values of the $\beta_{k,n}$s and later
the algorithm proceeds toward feasibility at the expense of some increase
of objective function values. However, even at those later sweeps
the objective function values of LinSup remain well below those of
the unsuperiorized application of the Cimmino feasibility-seeking
algorithm (Algorithm 2).

The slow increase of objective function values observed for the unsuperiorized
application of the Cimmino feasibility-seeking algorithm seems intriguing
because the feasibility-seeking algorithm is completely unaware of
the given objective function $\phi(x):=\left\langle c,x\right\rangle .$
But this is understood from the fact that the unsuperiorized algorithm
has an orbit of iterates in $R^{J}$ which, by proceeding in space
toward proximity minimizers, crosses the linear objective function's
level sets in a direction that either increases or decreases objective
function values. It would keep them constant only if the orbit was
confined to a single level set of $\phi$ which is not a probable
thing to happen. To clarify this we recorded in Figure 3 the values
of $\left\langle c,x\right\rangle $ and $\left\langle -c,x\right\rangle $
at the iterates $x^{k}$ produced by the Cimmino feasibility-seeking
algorithm (Algorithm 2).

\begin{figure}
\begin{raggedright} \includegraphics[scale=0.35]{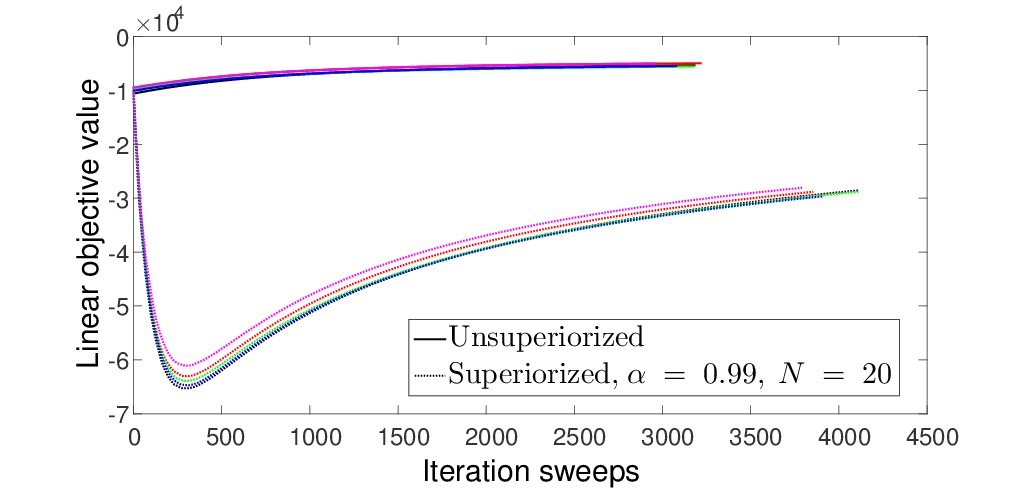}\caption{Linear objective function values plotted against iteration sweeps.
LinSup has reduced objective function values although the effect of
objective function reducing perturbations diminishes as iterations
proceed.\label{fig: lin-obj}}

\end{raggedright}
\end{figure}

\begin{figure}
\includegraphics[scale=0.35]{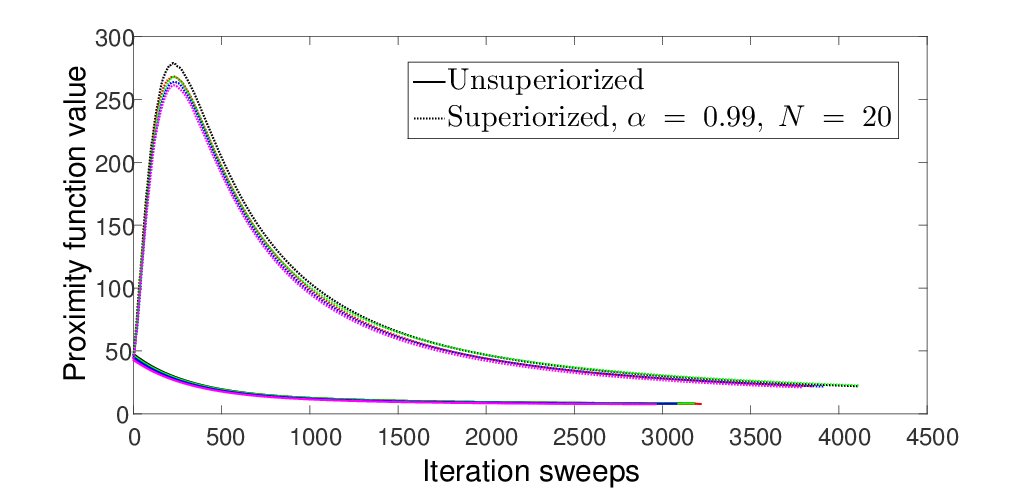}\caption{Proximity function values plotted against iteration sweeps. The unsuperiorized
feasibility-seeking only algorithm does a better job than LinSup here
which is understandable. LinSup's strive for feasibility comes at
the expense of some increase in objective function values, as seen
in Figure 1.\label{fig: prox}}
\end{figure}

\begin{figure}
\includegraphics[scale=0.35]{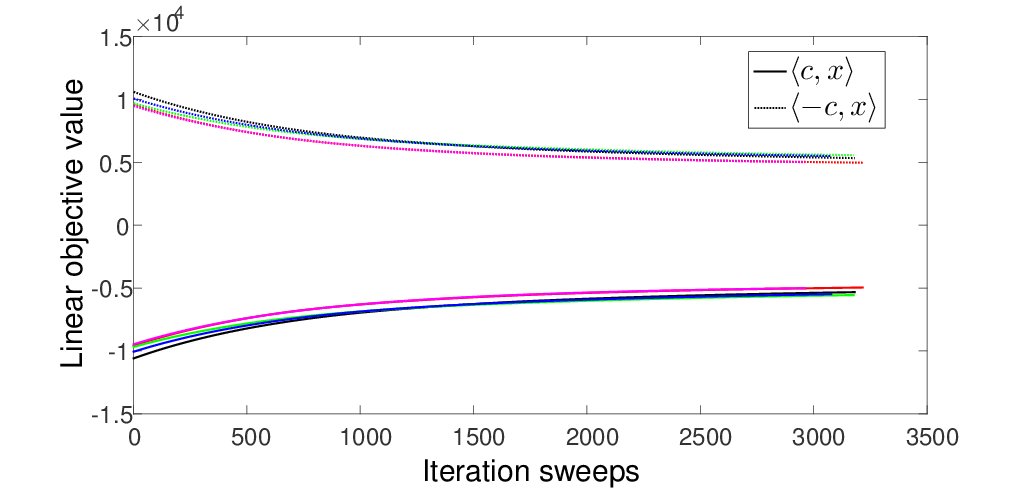}\caption{The fact that objective function values increase to some extent by
the unsuperiorized feasibility-seeking only algorithm observed in
Figure 1 is due to the relative situation of the linear objective
function's level sets with respect to where in space is the set of
proximity minimizers of the infeasible target set.\label{fig: c-c}}
\end{figure}

\section*{Concluding Comments}

We proposed a new approach to handle infeasible linear programs (LPs)
via the linear superiorization (LinSup) method. To this end we applied
the feasibility-seeking projection method of Cimmino to the original
linear infeasible constraints (without using additional variables).
This Cimmino method is guaranteed to converge to one of the points
that minimize a proximity function that measures the violation of
all constraints. We used the given linear objective function to superiorize
Cimmino's method to steer its iterates to proximity minimizers with
reduced objective function values. Further computational research
is needed to evaluate and compare the results of this new approach
to existing solution approaches to infeasible LPs.

\subsubsection*{Acknowledgments. \textmd{We thank Gabor Herman, Ming Jiang and Evgeni
Nurminski for reading a previous version of the paper and sending
us comments that helped improve it. This work was supported by Research
Grant No. 2013003 of the United States-Israel Binational Science Foundation
(BSF).}}

 \bibliographystyle{plain}
\bibliography{vvo-paper-bib-070416}

\end{document}